\documentstyle[11pt]{article}
\input amssym.def
\input amssym.tex

\newcount\skewfactor
\def\mathunderaccent#1#2 {\let\theaccent#1\skewfactor#2
\mathpalette\putaccentunder}
\def\putaccentunder#1#2{\oalign{$#1#2$\crcr\hidewidth
\vbox to.2ex{\hbox{$#1\skew\skewfactor\theaccent{}$}\vss}\hidewidth}}
\def\name{\mathunderaccent\tilde-3 }

\newcommand{\rest}{{\restriction}}

\newcommand{\dom}{{\rm dom}}

\newcommand{\forces}{\Vdash} 
\newcommand{\lesdot}{\mathrel{\mathord{<}\!\!\raise 0.8
pt\hbox{$\scriptstyle\circ$}}}  

\renewcommand{\a}{{\frak a}}
\renewcommand{\b}{{\frak b}}

\newcommand{\Proof}{\noindent {\sc Proof} \hspace{0.2in}} 
\newcommand{\pre}[2]{{}^{#1}{#2}}

\newcommand{\conc}{^\frown\!}

\newtheorem{theorem}{Theorem}

\newtheorem{main claim}[theorem]{Main Claim}
\newtheorem{lemma}[theorem]{Lemma} 
 
\newtheorem{corollary}[theorem]{Corollary} 
 
\newtheorem{definition}[theorem]{Definition}

\newtheorem{remark}[theorem]{Remark} 
\newtheorem{fact}[theorem]{Fact}

\makeatletter
\newcounter{abc}

\newcommand{\labelabc}{(\alph{abc})}
\newenvironment{abc}{\list{\labelabc}{\usecounter{abc}%
      \advance\@enumdepth \@ne 
	\def\makelabel##1{\hss\llap{##1}}}}{\endlist}

\newcommand{\pcf}{{\rm pcf}}

\def\cases#1{\left \{\,\vcenter {\normalbaselines \m@th \ialign
{$##\hfil $&\quad ##\hfil \crcr #1\crcr }}\right .}
\makeatother

\title{Embedding Cohen algebras using pcf theory}
\author{
 {\bf Saharon Shelah}\thanks{The research partially supported by ``The Israel
Science Foundation'' administered by of The Israel Academy of Sciences
and Humanities.  Publication 595.}\\
 Institute of Mathematics\\
 The Hebrew University of Jerusalem\\
 91904 Jerusalem, Israel\\
 and\\
 Department of Mathematics\\
 Rutgers University\\
 New Brunswick, NJ 08854, USA
 }
\date{done: July 1995\\
      printed: \today}

\begin{document}
\baselineskip13.14 truept

\maketitle

\begin{abstract} 
Using a theorem from pcf theory, we show that for any singular cardinal $\nu$,
the product of the Cohen forcing notions on $\kappa$, $\kappa<\nu$ adds a
generic for the Cohen forcing notion on $\nu^+$.
\end{abstract}

The following question (problem 5.1 in Miller's list \cite{Mi91}) is
attributed to Rene David and Sy Friedman:
\begin{quote} Does the product of the forcing notions
$\pre{{\aleph_n}>}2$ add a generic for the forcing
$\pre{{\aleph_{\omega+1}>}}2$? 
\end{quote}

We show here that the answer is yes in ZFC.   
Previously  Zapletal \cite{Za} has shown this result under the
assumption  $\square_{{\aleph_{\omega+1}}}$.

In fact, a similar theorem can be shown about other singular cardinals as
well.  The reader who is interested only in the original problem should read
$\aleph_{\omega+1}$ for $\lambda$, $\aleph_{\omega}$ for $\mu$ and $\{
\aleph_n: n\in(1,\omega)\}$ for $\a$.

\begin{definition}
\begin{enumerate}
\item 
Let $\a$ be a set of regular cardinals. $\prod\a$ is the set of all
functions $f$ with domain $\a$ satisfying $f(\kappa)\in \kappa$ for
all $\kappa \in \a$.
\item  A set $\b \subseteq \a $ is bounded if $\sup \b < \sup\a$.
$\b$ is cobounded if $\a \setminus \b $ is bounded. 
\item    If $J$ is an ideal on $\a$, $f,g\in \prod \a$, then $f<_J g$
means $\{\kappa \in \a: f(\kappa)\not< g(\kappa) \} \in J$. 
We write $\prod \a/J$ for the partial (quasi)order $(\prod \a,{ <_J})$. 
\item $\lambda = tcf(\prod \a/J)$ ($\lambda $  is the true cofinality 
of $\prod \a/J$) means that there is a strictly increasing cofinal
sequence of functions in the partial order $(\prod \a, <_J)$. 
\item $\pcf(\a) = \{ \lambda: \exists J \,\, \lambda =
tcf(\prod\a/J)\}$. %  (We may restrict $J$ to be a maximal ideal.)
\end{enumerate}
\end{definition}

We will use the following theorem from pcf theory:

\begin{lemma}\label{pcf}
Let $\mu$ be a singular cardinal.  Then there is a 
set $\a$  of regular 
cardinals below $\mu$, $|\a|= cf(\mu) < \min \a$ 
and $\mu^+\in \pcf(\a)$. 

Moreover, we can even have 
$tcf(\prod \a/J^{bd})=\mu^+$, where $J^{bd}$ is the ideal 
of all bounded subset of $\a$. 
\end{lemma}

\Proof  See \cite[theorem 1.5]{Sh:355}.

\begin{theorem}\label{maintheorem}
Let $\a$ be a set of regular cardinals, $\mu=\sup \a\notin \a$,
$2^{<\lambda } = 2^\mu$, $\lambda > \mu$, 
$\lambda \in \pcf(\a)$, and moreover:  
\begin{itemize}
\item[$(*)$] There is an ideal $J$ on $\a$ containing all 
bounded sets such that 
$\lambda =tcf(\prod \a/J)$.  
\end{itemize}
Then the forcing notion $\prod _{\kappa\in\a} \pre{\kappa >}2$
adds a generic for $\pre{\lambda >}2$. 
\end{theorem}

\begin{corollary}
If $\nu$ is a singular cardinal, and $P$ is the product of the
forcing notions $\pre{\kappa >}2$ for $\kappa < \nu$, then $P$ adds
a generic for $\pre{\nu^+>}2$.
\end{corollary}
\Proof  By lemma \ref{pcf} and theorem \ref{maintheorem}

\begin{remark}
\begin{enumerate}
\item The condition $(*)$ in the theorem is equivalent to:
\begin{itemize}
\item[$(**)$] For all bounded sets $\b \subset \a$ 
  we have $\lambda \in \pcf(\a \setminus \b)$.
\end{itemize}
\item  Clearly the assumption 
$2^{<\lambda } = 2^\mu$ is necessary, because otherwise the 
forcing notion 
 $\prod _{\kappa\in\a} \pre{\kappa >}2$ would be too small to 
add a generic for $\pre{\lambda >}2$.
\end{enumerate}
\end{remark}

\def\ppp{\prod \a } 
\def\liff{\Leftrightarrow}
\def\nr{\name \rho}

\noindent {\bf Proof of the theorem:}
By our assumption we have some ideal $J$ containing 
all bounded sets such that $tcf(\ppp / J) = \lambda$.  

We will write $\forall ^{J} \kappa \in \a\,\, \varphi(\kappa)$ for
$$ \{ \kappa \in \a : \lnot \varphi(\kappa) \} \in J$$

So we have a 
 sequence 
 $\langle f_\alpha :
\alpha < {\lambda} \rangle $ such that 
\begin{abc}
\item $f_\alpha \in \ppp $
\item If $\alpha < {\beta}$, then $\forall^{J}  \kappa \in  \a \,\, 
f_\alpha( \kappa ) < f_{\beta}( \kappa )$
\item $\forall f\in \ppp \, \, \exists \alpha \, \,      
\forall^{J}  \kappa \in  \a \,\, 
f( \kappa ) < f_\alpha( \kappa )$.
\end{abc}

The next lemma shows that if we allow these fucntions to be defined
only almost everywhere, then we can additionally assume that 
in each block of length $\mu$  these functions
have  disjoint graphs:

\medskip

\begin{lemma}\label{pcf2}
Assume that $\a$, $\lambda$, $\mu$ are as above.  Then 
there is a 
sequence $\langle g_\alpha: \alpha < {\lambda} \rangle$ such
that 
\begin{abc}
\item 	\label{cob}
$\dom(g_\alpha) \subseteq \a $ cobounded
(so in particular $\forall^J \kappa\in \a: \kappa\in  
\dom(g_\alpha(\kappa))$.  
\item If $\alpha  < {\beta} $, then $\forall^{J}  \kappa \in  \a  \,\, 
	g_\alpha ( \kappa ) < g_{\beta} ( \kappa )$\label{inc}
\item  $\forall f\in \ppp \, \, \exists \alpha  \, \,   
\forall^{J}  \kappa \in  \a \,\,  f( \kappa ) < g_\alpha ( \kappa )$.
Moreover, we may   
choose $\alpha$ to be divisible by ${\mu}$. \label{cof} 
\item If $\alpha < {\beta} < \alpha  + {\mu} $, then
$\forall  \kappa \in \dom(g_\alpha)\cap \dom(g_{\beta}):
g_\alpha(\kappa) < g_{\beta}(\kappa)$. 
 \label{disjoint}
\end{abc}
\end{lemma}

\Proof  Let   
$\langle f_{\alpha} : {\alpha} < {\lambda}
\rangle$ be as above.  Now define $\langle 
g_\alpha            : \alpha < {\lambda} \rangle $ by induction 
as follows: \\
  If $\alpha = {\mu} \cdot {\zeta} $,  then let $g_{\alpha}$ 
be any function that satisfies $g_\beta <_J g_\alpha$ for all 
$\beta < \alpha$, and also $f_\alpha <_J g_\alpha$.   Such a function
can be found because set of functions of size $< \lambda$ can be
$<_J$-bounded by some $f_{\beta}$. 
\\
If  $\alpha = {\mu} \cdot {\zeta}+i $, $0<i<\mu$, 
  then let 
$$ g_\alpha( \kappa ) =
 \cases{ g_{\mu\cdot \zeta}(\kappa) + i  & if  $i <{\kappa }$\cr
	\hbox{undefined} & otherwise \cr}$$
It is easy to see that (a)--(d) are satisfied.

\begin{definition}
\begin{enumerate}
\item Let $P_ \kappa $ be the set  $\pre{\kappa >}2$, partially ordered by
inclusion (= sequence extension). Let $P = \prod_{ \kappa \in  \a }
P_\kappa $. 
[We will show that $P$ adds a generic for $\pre{{\lambda} > }2$]
\item Assume that $\langle g_{\alpha}: {\alpha} < {\lambda}
\rangle$ is as in lemma~\ref{pcf2}. 
\item Let  $H : \pre {\mu}2 \to 
      \pre{\lambda >}2  $ be 
 onto. 
 \item For $\kappa\in \a$, 
let $\name\eta_\kappa$ be the $P_ \kappa $-name for the  generic function
from $\kappa$ to $2$.  Define a $P$-name 
 of a 
function $\name h :{\lambda}  \to 2$
by 
$$  \name h(\alpha) = 
\cases{ 0 & if $
\forall^ {J}  
\kappa  \in  \a  \, \, \name \eta_ \kappa (g_\alpha( \kappa )) = 0$\cr
1& otherwise } $$
\item For ${\xi} < {\lambda}$ let 
 $\nr_{\xi}$ be a $P$-name for the element  of $ \pre {\mu}2$
that satisfies 
$\nr_{\xi} \simeq \name h \rest [{\mu} \cdot {\xi},
{\mu} \cdot ({\xi} +1 ))$, i.e., 
$$
i < {\mu}\ \Rightarrow \ 
\forces_{P}\, 
\nr_{\xi}(i)=  
\name h({{\mu} \cdot {\xi} + i}).$$ 

 Define $\nr \in \pre
{{\lambda}}2$ by 
$$ \nr = 
H(\nr_0) \conc H(\nr_1) \conc \cdots \conc H(\nr_{\xi}) \conc
\cdots $$
\end{enumerate}
\end{definition}

\begin{main claim}  $\nr$ is generic for $\pre{{\lambda}>}2$. 
\end{main claim}

\begin{definition}
For $\alpha < {\lambda}$ let $P^{(\alpha)}$ be the set of
all conditions $p$ satisfying $\forall^{J}  \kappa :  \dom(p_ \kappa )=
g_\alpha( \kappa )$. 
\end{definition}

\begin{remark}
$\bigcup_{{\zeta}  < {\lambda}} P^{({\mu} \cdot
{\zeta} )}$ is dense in $P$. 
\end{remark}
\Proof By lemma~\ref{pcf2}\ref{cof}. 

\begin{fact}\label{extend}
Let $\alpha = {\mu} \cdot {\zeta}$, 
 $p \in P^{(\alpha)}$, $\sigma \in \pre {\mu} 2$.  Then
there is a condition $q\in P^{(\alpha+{\mu})}$, $q \ge p$ and 
$$ \forall j < {\mu} \,\, \forall^{J}  \kappa  \,\,
	q_ \kappa (g_{\alpha+j}( \kappa )) = \sigma(j)$$
\end{fact}
\Proof   Let $p= (p_ \kappa :  \kappa  \in \a)$. 
There is a set $\b\in J$ such that: For all 
  $\kappa\in \a \setminus \b$ 
   we have $\dom(p_\kappa) =
g_{\alpha}(\kappa)$.   Define  $q\in P^{(\alpha+{\mu})}$,
 $q = (q_ \kappa :  \kappa  \in \a)$  as follows: 
$$q_ \kappa ({\gamma}) = \cases{
	p_ \kappa ({\gamma})   & if ${\gamma} \in \dom(p_ \kappa )$\cr
	\sigma(j)               & if ${\gamma} = g_{\alpha+j}( \kappa )$, 
			 $\kappa\in \a \setminus \b$\cr
	0               & otherwise \cr}
$$
We have to explain why $q$ is well-defined:   First note that the
first and the second case are mutually exclusive.   Indeed, if
${\gamma} = g_{\alpha+j}( \kappa )$, then ${\gamma} >
 g_\alpha(\kappa)$, whereas $\kappa \notin \b$ implies that 
$\dom(p_\kappa)=g_\alpha(\kappa)$, so $\gamma\notin\dom(p_\kappa)$. 
\\
Next, by 
the property~\ref{disjoint} from 
lemma~\ref{pcf2} there is no contradiction between various
instances of the second case.

Hence we get that for all $j < {\mu}$, whenever $\kappa \in \a 
\setminus \b$, 
and ${\kappa } > j$, then $q_ \kappa (g_{\alpha+j}( \kappa )) =
\sigma(j)$.   Since $J$ contains all bounded sets, this means that  
$\forall^J \kappa : q_ \kappa (g_{\alpha+j}( \kappa )) =
\sigma(j)$.   

\begin{remark}\label{extendrem}
Assume that $\alpha = {\mu} \cdot {\zeta}$, and $p$, $q$,
${\sigma}$ are as above.  Then  
$q \forces \nr_{\zeta} = {\sigma}$. 
\end{remark}
\medskip

\noindent {\bf Proof of the main claim:}
Let $p\in P$, and $D \subseteq \pre{{\lambda}>} 2 $ be a
dense open set.  We may assume that for some $\alpha^*<
{\lambda} $, ${\zeta}^* <
{\lambda}$ we have $\alpha^* = {\mu} 
\cdot {\zeta}^*$ and $p \in P^{(\alpha^*)}$, i.e., for some
 $\b\in J$:  
$ \forall  \kappa \notin \b :
\dom(p_ \kappa ) = g_{\alpha^*}( \kappa ) 
$

So $p$ decides the values of $h \rest \alpha^*$, and hence also the
values of $\nr_{{\zeta} }$ for ${\zeta} < {\zeta}^*$. Specifically,
for each ${\zeta} < {\zeta}^*$ we can define ${\sigma}_{\zeta} \in
 \pre {\mu} 2$ by 
$$ {\sigma}_{\zeta}(i) = \cases { 0 & if $\forall ^{J}  \kappa \,\,
p_ \kappa (g_{{\mu} \cdot {\zeta} + i}( \kappa )) = 0 $\cr
1 & otherwise \cr}$$
(Note that for all ${\zeta} < {\zeta}^*$, for all $i <
{\mu}$, for almost all $ \kappa $ the value of 
$p_ \kappa (g_{{\mu} \cdot {\zeta} + i}( \kappa ))$ is defined.)

Clearly $p \forces \nr_{\zeta} = {\sigma}_{\zeta}$. 
Since $D$ is dense and $H$ is onto, we can 
now find ${\sigma}_{{\zeta}^*}\in \pre {\mu} 2$ such
that 
$$ H({\sigma}_0) \conc \cdots \conc H({\sigma}_{\zeta}^*) 
		\in D$$
Using \ref{extend} and \ref{extendrem}, we can now find $q \ge p$ such
that $q \forces \rho_{{\zeta}^*} = {\sigma}_{\zeta^*}$. 

Hence $q \forces \nr\in D$, and we are 
%                         ----------------
			       done.
%                         ----------------

%\bibliographystyle{lit-unsrt}
%\bibliography{listb,listx}

\end{document}